\documentclass[a4paper,12pt]{amsart}
\usepackage[letterpaper, margin=1.2in]{geometry}
\usepackage{latexsym}
\usepackage{amssymb}
\usepackage{amsmath}
\usepackage{amsfonts}
\usepackage[table]{xcolor}
\usepackage{pgfplots}
\usepackage{color}
\usepackage{multirow}
\usepackage{graphicx}
\usepackage{txfonts,pxfonts} 
\usepackage{tikz}
\usetikzlibrary{calc,arrows}
\usetikzlibrary{calc}
\usepackage{verbatim}
\everymath{\displaystyle}
\newtheorem{thm}{Theorem}[section]

\newtheorem{lem}[thm]{Lemma}

\newtheorem{cor}[thm]{Corollary}

\newtheorem{prop}[thm]{Proposition}

\theoremstyle{definition}

\newtheorem{rem}{Remark}

\newtheorem*{example}{Example}

\def\al{\alpha}
\def\la{\lambda}

\def\th{\theta}
\def\md#1{\ \mbox{\rm(mod }{#1})}

\newcounter{cs}
\stepcounter{cs}
\newcommand{\casos}{\begin{itemize}}
\newcommand{\fcasos}{\end{itemize}\setcounter{cs}{1}}

\newfont{\tit}{cmr12 scaled \magstep3}
 \title{On power integral bases of certain pure
 	number fields defined by $X^{60}-m$ }
 
 \author{{Lhoussain El FADIL $^{1}$, Omar KCHIT$^{ 2}$, and  Hanan CHOULLI$^{ 3}$}
    $^{1}$Department of Mathematics, Faculty of Sciences Dhar El Mahraz, Sidi Mohamed ben Abdellah University, Fes, Morocco,\\ ORCID iD:  https://orcid.org/0000-0003-4175-8064\\
 	$^{2}$Department of Mathematics, Faculty of Sciences Dhar El Mahraz, Sidi Mohamed ben Abdellah University, Fes, Morocco,\\ ORCID iD:  https://orcid.org/0000-0002-0844-5034\\
 	$^{3}$Department of Mathematics, Faculty of Sciences Dhar El Mahraz, Sidi Mohamed ben Abdellah University, Fes, Morocco, \\  ORCID iD: https://orcid.org/0000-0002-5045-9069\\
 	\\ 
e-mail: lhoussain.elfadil@usmba.ac.ma 
  }
 		
 
\begin{document}

 	\begin{abstract}  Let $K$ be a pure number field generated by a complex root of a monic irreducible polynomial $F(x)=x^{60}-m\in \mathbb{Z}[x]$, with $m\neq \pm1$  a square free integer. In this paper, we study the monogeneity of $K$. We prove that if $m\not\equiv 1\md{4}$, $m\not\equiv \mp 1 \md{9} $ and $\overline{m}\not\in\{\mp 1,\mp 7\} \md{25}$, then  $K$ is monogenic. But if  $m\equiv 1\md{4}$, $m\equiv \mp1 \md{9}$, or $m\equiv \mp 1\md{25}$, then  $K$ is not  monogenic. Our results are illustrated by examples.
 	   	\end{abstract}
 \maketitle
 	 	
 	\section{\bf Introduction}\label{sec:1}
Let $K=\mathbb{Q}(\alpha)$ be a number field generated by  a complex root $\alpha$ of a monic irreducible polynomial $F(x)\in \mathbb{Z}[x]$ and $\mathbb{Z}_K$ its ring of integers. It is well known that the ring $\mathbb{Z}_K$ is a free $\mathbb{Z}$-module of rank $n=[K:\mathbb{Q}]$, and so the Abelian group $\mathbb{Z}_K/\mathbb{Z}[\alpha]$ is finite. Its cardinal order  is called the index of $\mathbb{Z}[\alpha]$ and denoted $(\mathbb{Z}_K:\mathbb{Z}[\alpha])$. 
If for some generator $\th\in\mathbb{Z}_K$ of $K$, we have $(\mathbb{Z}_K:\mathbb{Z}[\th])=1$, then the ring $\mathbb{Z}_K$ is said to have a power integral basis $(1,\theta,\dots,\theta^{n-1})$. In such a case, the field $K$ is said to be monogenic and not monogenic otherwise. The problem of testing the monogeneity of number fields and constructing power integral bases
have been intensively studied these last fourth decades, mainly by Ga\'al, Gy\"ory, Nakahara,
Pohst and their research teams (see for instance \cite{AN, 2a, G, G19, P}).
It is called a problem of Hasse to give an arithmetic characterization of those number fields which have a power integral basis \cite{2a, F4, Ha, He, MNS, P}.
{ In \cite{E07},  El Fadil  gave conditions for the existence of power integral bases of pure cubic fields in terms of the index form equation. In \cite{F4}, Funakura, calculated integral bases of pure quartic fields and studied their monogeneity.  In \cite{GR4},  Ga\'al and  Remete, calculated the elements of index $1$ of pure quartic fields generated by $m^{\frac{1}{4}}$ for {$1< m <10^7$}  and $m\equiv 2,3 \md{4}$.  In \cite{AN6}, Ahmad, Nakahara, and Husnine proved  that  if $m\equiv 2,3 \md{4}$ and  $m\not\equiv \mp1\md{9}$, then the sextic number field generated by $m^{\frac{1}{6}}$ is monogenic.
	They also showed in \cite{AN},    that if $m\equiv 1 \md{4}$ and $m\not\equiv \mp1\md{9}$, then the sextic number field generated by $m^{\frac{1}{6}}$ is not monogenic.  In \cite{E6}, based on prime ideal factorization, El Fadil showed that  if $m\equiv 1 \md{4}$ or $m\not\equiv 1\md{9}$, then the sextic number field generated by $m^{\frac{1}{6}}$ is not monogenic.
	{Hameed and Nakahara \cite{HN8}, proved that if $m\equiv 1\md{16}$, then the octic number field generated by $m^{1/8}$ is not monogenic, but if $m\equiv 2,3 \md{4}$, then it is monogenic.}   In \cite{GR17}, by applying the  explicit form of the index equation, Ga\'al  and  Remete obtained deep new results  and they gave a complete answer  to the problem   of monogeneity of  number  fields generated  by $m^{\frac{1}{n}}$, where $3\leq n\leq 9$.  While Ga\'al's and  Remete's techniques are based on the index calculation,  Nakahara's methods are based on the existence of power relative integral bases of some special sub-fields.  In \cite{ E6, E6s, E12, E18, E24, E36},  El Fadil et al.  used   Newton polygon techniques to study the monogeneity of the pure number fields of degrees $6,  12, 18, 24$, and $36$.
	In this paper,
Our purpose is for a square free integer $m\neq \pm 1$ and $F(x)=x^{60}-m$ is an irreducible polynomail over $\mathbb{Q}$, to study the monogeneity of  the number field $K=\mathbb{Q}(\alpha)$ generated by a complex root  $\alpha$ of  $F(x)$.
Our method applying the  Newton polygon techniques and the explicit prime ideal factorization.
Recall that some authors like Khanduja and Gassert  introduced a new topic of monogeneity of polynomials as follows (see \cite{Kh, G}): For a  monic irreducible polynomial $F(x)$ and $K$ the number field generated by a root of $F(x)$,  is said to be monogenic if the ring $\mathbb{Z}_K=\mathbb{Z}[\alpha]$, and not monogenic otherwise. Unfortunately, this notion of monogeneity is not noting other the integral closedness of $\mathbb{Z}[\alpha]$. Even if the problem of the integral closedness of $\mathbb{Z}[\alpha]$ gives a partial answer to the problem of monogeneity of $K$, the problem of monogeneity is more hard to solve. Unfortunately this notion of monogeneity  does not coincide with the well known notion of monogeneity of number fields,  treated by Gaal, Nakahara, Pohst and their collaborators. In fact let us consider  $F(x)=x^3-9$ and $K$ the number field generated by a complex root  $\alpha$  of $F(x)$,  since $3$ divides the index $(\mathbb{Z}_K:\mathbb{Z}[\alpha])$, we conclude that  $\mathbb{Z}[\alpha]$ is not integrally closed even if  $K$ is monogenic and  $\theta=\frac{\alpha^2}{3}$ generates a power integral closure of $K$ (For more details on the confusion between the two different notions of monogeneity, we refer to El Fadil's preprint “A note on monogeneity of pure number fields” on the link https://arxiv.org/abs/2106.00004). 

\section{Main Results}
Let $K$ be a pure number field defined by a complex root $\alpha$ of a monic irreducible polynomial $F(x)=x^{60}-m$, with $m\neq \pm 1$  a square free integer.
The main goal of this section is  to study the monogeneity of pure number fields of degree $60$. Theorem \ref{thm1} gives a necessary and sufficient condition on the integral closure of $ \mathbb{Z}[\alpha]$. This  theorem covers  \cite[Theorem 1.1]{G} in the context of pure number fields of degree $60$. But on the contrary, \cite[Theorem 1.1]{G} does not cover our Theorem \ref{thm1}. In fact, our theorem gives necessary and sufficient conditions on the integral closure of $\mathbb{Z}[\alpha] $, unlike Gassert's results, which gives just one meaning and requires more details to reach our result. Theorem \ref{thm2}  gives a partial converse on Theorem \ref{thm1}. In fact it gives a full converse for the number fields defined by a monic irreducible polynomial $F(x)=x^{60}-m$, with $m$ a square free integer, except for the cases $m\in \{+7, -7\} \md{25}$. Finally Theorem \ref{cor} gives a partial answer for  the cases when $m$ is not necessarily a square free integer.
\begin{thm}\label{thm1}
	The ring $\mathbb{Z}[\alpha]$ is the ring of integers of $K$ if and only if $m\not\equiv 1\md {4}$, $m\not\equiv \mp 1\md {9}$ and $m\not\equiv \mp1,\mp 7\md {25}$.\\
	In particular, if $m\not\equiv 1\md {4}$, $m\not\equiv \mp 1\md {9}$ and $m\not\equiv \mp1,\mp 7\md {25}$, then $K$ is monogenic.
\end{thm}
Remark that based on Theorem \ref{thm1}, if $m\equiv 1\md{4}$ or $m\equiv \mp1\md{9}$ or $m\equiv \mp1,\mp 7\md {25}$ then $\mathbb{Z}[\alpha]$ is not the ring of integers of $K$. But in this case, Theorem \ref{thm1} can not decide on monogeneity of $K$. The following theorem gives a partial answer, the exception is when $m\equiv \mp 7\md {25}$.
\begin{thm}\label{thm2}
	If one of the following statements holds,
	\begin{enumerate}
		\item
		$m\equiv 1 \md{4}$,
		\item
		$m\equiv \mp 1\md{9}$, 
		\item 
		$m\equiv \mp 1\md{25}$,
	\end{enumerate}
	then $K$ is not monogenic.
\end{thm}
\begin{thm}\label{cor}
	Let $K$ be a pure number field defined by a complex root $\alpha$ of a monic irreducible polynomial $F(x)=x^{60}-a^u$, with $a\neq \mp 1$  a square free integer and $u$ a positive integer which is coprime to $30$. Then 
	\begin{enumerate}
		\item
		If $a\not\equiv 1\md{4}$, $a\not\equiv \mp1\md{9}$, and $a\not\equiv \mp1,\mp 7\md {25}$,  then  $K$ is monogenic.
		\item 
		If $a\equiv 1\md{4}$ or  $a\equiv \mp 1\md{9}$ or $a\equiv \mp1\md{25}$, 
		then $K$ is not monogenic.
	\end{enumerate}
\end{thm}
\section{Preliminaries}
In order to prove our main Theorems, we recall some fundamental facts about Newton polygon techniques applied on prime ideal factorization and calculation of index. 
Let  $\overline{F(x)}=\prod_{i=1}^r \overline{\phi_i(x)}^{l_i}$ in  $\mathbb{F}_p[x]$ be the factorization of $\overline{F(x)}$ into powers of monic irreducible coprime polynomials of $\mathbb{F}_p[x]$. Recall that  a well-known theorem of  Dedekind says that: 
\begin{thm}$($\cite[ Chapter I, Proposition 8.3]{Neu}$)$\\
	If $p$ does not divide the index $(\mathbb{Z}_K:\mathbb{Z}[\al])$, then
	$$p\mathbb{Z}_K=\prod_{i=1}^r \mathfrak{p}_i^{l_i}, \mbox{ where every } \mathfrak{p}_i=p\mathbb{Z}_K+\phi_i(\al)\mathbb{Z}_K,$$
	and the residue degree of $\mathfrak{p}_i$ is $f(\mathfrak{p}_i)={\mbox{deg}}(\phi_i)$.
\end{thm}
In order to apply this  theorem in an effective way one needs a
criterion to test  whether  $p$ divides  or not the index $(\mathbb{Z}_K:\mathbb{Z}[\al])$.
In $1878$, Dedekind  gave a criterion  to test whether  $p$ divides or not  $(\mathbb{Z}_K: \mathbb{Z}[\al])$.
{
	\begin{thm}\label{Ded}$($Dedekind's Criterion \cite[Theorem 6.1.4]{Co} and \cite{R}$)$\\
		For a number field $K$ generated by $\al$ a complex root of a monic irreducible  polynomial $F(x)\in \mathbb{Z}[x]$ and a rational prime integer $p$, let $\overline{F(x)}=\prod_{i=1}^r\overline{\phi_i(x)}^{l_i}\md{p}$  be the factorization of   $\overline{F(x)}$ in $\mathbb{F}_p[x]$, where the polynomials $\phi_1(x),\dots,\phi_r(x)$ are monic polynomials in $\mathbb{Z}[x]$, their reductions are coprime polynomials irreducible over $\mathbb{F}_p$. If we set
		$M(x)=\cfrac{F(x)-\prod_{i=1}^r{\phi_i}(x)^{l_i}}{p}$, then $M(x)\in \mathbb{Z}[x]$ and the following statements are equivalent:
		\begin{enumerate}
			\item[1.]
			$p$ does not divide the index $(\mathbb{Z}_K:\mathbb{Z}[\al])$.
			\item[2.]
			For every $i=1,\dots,r$, either $l_i=1$ or $l_i\geq 2$ and $\overline{\phi_i(x)}$ does not divide $\overline{M}(x)$ in $\mathbb{F}_p[x]$.
		\end{enumerate}
\end{thm} }
When Dedekind's criterion fails, that is,  $p$ divides the index $(\mathbb{Z}_K:\mathbb{Z}[\th])$ for every primitive element $\th\in\mathbb{Z}_K$,{  then} for such primes and number fields, it is not possible to obtain the prime ideal factorization of $p\mathbb{Z}_K$ by {Dedekind's theorem}.
In 1928, Ore developed   an alternative approach
for obtaining the index $(\mathbb{Z}_K:\mathbb{Z}[\alpha])$, the
absolute discriminant $d_K$, and the prime ideal factorization of the rational primes in
a number field $K$ by using Newton polygons (see for instance {\cite{EMN, MN, O}}). 
Now we recall  some fundamental facts  on Newton polygons, for more details, we refer to \cite{El,GMN}. For any prime integer $p$ and for any monic polynomial 
$\phi\in
\mathbb{Z}[x]$  whose reduction is irreducible  in
$\mathbb{F}_p[x]$, let $\mathbb{F}_{\phi}$ be the finite field $\mathbb{F}_p[x]/(\overline{\phi})$. For any monic polynomial $F(x)\in \mathbb{Z}[x]$, upon to the Euclidean division by successive powers of $\phi$, we expand $F(x)$ as $F(x)=a_0(x)+a_1(x)\phi(x)+\cdots+a_l(x)\phi(x)^l$, called the $\phi$-expansion of $F(x),~($for every $i,~\deg(a_i(x))<\deg(\phi))$. To any coefficient $a_i(x)$ we attach $u_i=\nu_p(a_i(x))\in \mathbb{Z}\cup \{\infty\}$. The $\phi$-Newton polygon of $F(x)$ with respect to $p$, is the lower boundary convex envelope of the set of points $\{(i,u_i),~ a_i(x)\neq 0\}$ in the Euclidean plane, which we denote by $N_{\phi}(F)$. Geometrically, the $\phi$-Newton polygon of $F$, is the process of joining the obtained edges  $S_1,\dots,S_t$ ordered by increasing slopes, which  can be expressed as $N_{\phi}(F)=S_1+\dots + S_t$.  The principal $\phi$-Newton polygon of ${F}$},
denoted $N_{\phi}^-(F)$, is the part of the  polygon $N_{\phi}(F)$, which is  determined by joining all sides of $N_{\phi}(F)$ of negative  slopes.
For every side $S$ of $N_{\phi}^-(F)$, { the length of $S$, denoted $l(S)$, is the length of its projection to the $x$-axis and  its height, denoted $h(S)$, is the length of its projection to the $y$-axis}. {Let $d=$gcd$(l(S), h(S))$ be the  degree of $S$.
For every side $S$ of {$N_{\phi}^-(F)$}, with initial point $(s, u_s)$, length $l$, and for every 
abscissa $i=0, \dots,l$, we attach   the following
{{\it residue coefficient} $c_i\in\mathbb{F}_{\phi}$ as follows:
	$$c_{i}=
	\left
	\{\begin{array}{ll} 0,& \mbox{ if } (s+i,{\it u_{s+i}}) \mbox{ lies strictly
		above } S,\\
	\left(\dfrac{a_{s+i}(x)}{p^{{\it u_{s+i}}}}\right)
	\,\,
	\mod{(p,\phi(x))},&\mbox{ if }(s+i,{\it u_{s+i}}) \mbox{ lies on }S.
	\end{array}
	\right.$$
	where $(p,\phi(x))$ is the maximal ideal of $\mathbb{Z}[x]$ generated by $p$ and $\phi$.
	Let $\lambda=-h/e$ be the slope of $S$, where  $h$ and $e$ are two positive coprime integers. Then  $d=l/e$ is the degree of $S$. Since 
	the points  with integer coordinates lying{ on} $S$ are exactly $${(s,u_s),(s+e,u_{s}-h),\cdots, (s+de,u_{s}-dh)},$$ if $i$ is not a multiple of $e$, then 
	$(s+i, u_{s+i})$ does not lie on $S$, and so $c_i=0$. Let
	{$$R_{\lambda}(F)(y)=t_dy^d+t_{d-1}y^{d-1}+\cdots+t_{1}y+t_{0}\in\mathbb{F}_{\phi}[y],$$}} called  
the residual polynomial of $F(x)$ associated to the side $S$, where for every $i=0,\dots,d$,  $t_i=c_{ie}$.\\
\smallskip

Remark that as $(s, u_{s})$ and $(s+l, u_{s+l})$  lie on $S$, then deg$(R_{\lambda}(F)(y))=d$ and $t_0\neq 0$.\\
\smallskip

Let $N_{\phi}^-(F)=S_1+\dots + S_t$ be the principal $\phi$-Newton polygon of $F$ with respect to $p$.\\
We say that $F$ is a $\phi$-regular polynomial with respect to $p$, if  $R_{\lambda_i}(F)(y)$ is square free in $\mathbb{F}_{\phi}[y]$ for every  $i=1,\dots,t$. 
The polynomial $F$ is said to be  $p$-regular  if $\overline{F(x)}=\prod_{i=1}^r\overline{\phi_i(x)}^{l_i}$ for some monic polynomials $\phi_1(x),\dots,\phi_r(x)$ of $\mathbb{Z}[x]$ such that $\overline{\phi_1(x)},\dots,\overline{\phi_r(x)}$ are irreducible coprime polynomials over $\mathbb{F}_p$ and $F$ is  a $\phi_i$-regular polynomial with respect to $p$ for every $i=1,\dots,r$.
\smallskip

The  theorem of Ore plays  a  key role for proving our main Theorems.\\
Let $\phi(x)\in\mathbb{Z}[x]$ be a monic polynomial, with $\overline{\phi(x)}$ is irreducible in $\mathbb{F}_p[x]$. As defined in \cite[Def. 1.3]{EMN},  the $\phi$-index of $F(x)$, denoted by $ind_{\phi}(F)$, is  deg$(\phi)$ times the number of points with natural integer coordinates that lie below or on the polygon $N_{\phi}^-(F)$, strictly above the horizontal axis,{ and strictly beyond the vertical axis} (see $Figure\ 1$).\\
\begin{figure}[htbp] 
	\centering
	\begin{tikzpicture}[x=1cm,y=0.5cm]
	\draw[latex-latex] (0,6) -- (0,0) -- (10,0) ;
	\draw[thick] (0,0) -- (-0.5,0);
	\draw[thick] (0,0) -- (0,-0.5);
	\node at (0,0) [below left,blue]{\footnotesize $0$};
	\draw[thick] plot coordinates{(0,5) (1,3) (5,1) (9,0)};
	\draw[thick, only marks, mark=x] plot coordinates{(1,1) (1,2) (1,3) (2,1)(2,2)     (3,1)  (3,2)  (4,1)(5,1)  };
	\node at (0.5,4.2) [above  ,blue]{\footnotesize $S_{1}$};
	\node at (3,2.2) [above   ,blue]{\footnotesize $S_{2}$};
	\node at (7,0.5) [above   ,blue]{\footnotesize $S_{3}$};
	\end{tikzpicture}
	\caption{The principal  $\phi$-Newton polygon of $F(x)$.}
\end{figure}
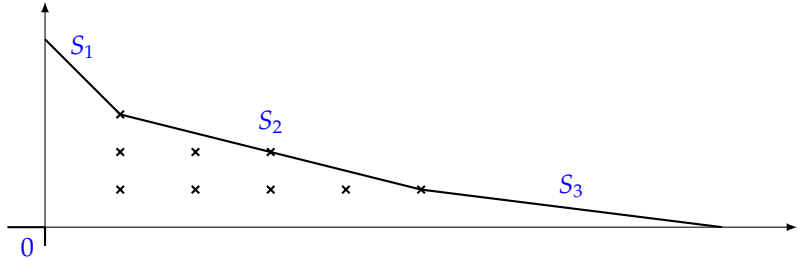
In the example of $Figure\ 1$, $ind_\phi(F)=9\times$deg$(\phi)$.\\
\smallskip

 Now assume that $\overline{F(x)}=\prod_{i=1}^r\overline{\phi_i}^{l_i}$ is the factorization of $\overline{F(x)}$ in $\mathbb{F}_p[x]$, where every $\phi_i\in\mathbb{Z}[x]$ is monic polynomial, with $\overline{\phi_i(x)}$ is irreducible in $\mathbb{F}_p[x]$, $\overline{\phi_i(x)}$ and $\overline{\phi_j(x)}$ are coprime when $i\neq j$ and $i, j=1,\dots,r$.
For every $i=1,\dots,r$, let  ${N_{\phi_i}^-(F)=S_{i1}+\dots+S_{it_i}}$ be the principal  $\phi_i$-Newton polygon of $F(x)$ with respect to $p$. For every ${j=1,\dots, t_i}$,  let $R_{\lambda_{ij}}(F)(y)=\prod_{s=1}^{s_{ij}}\psi_{ijs}^{a_{ijs}}(y)$ be the factorization of $R_{\lambda_{ij}}(F)(y)$ in $\mathbb{F}_{\phi_i}[y]$. 
Then we have the following index theorem of Ore (see \cite[Theorem 1.7 and Theorem 1.9]{EMN}, \cite[Theorem 3.9]{El}, and{\cite[pp: 323--325]{MN}}).
\begin{thm}\label{thm4}$($Theorem of Ore$)$ Under the above hypothesis,
	\begin{enumerate}
		\item 
		$$\nu_p((\mathbb{Z}_K:\mathbb{Z}[\alpha]))\geq\sum_{i=1}^{r}ind_{\phi_i}(F).$$  
		The equality holds if $F(x) \text{ is }p$-regular.
		\item 
		If  $F(x) \text{ is }p$-regular, then
		$$p\mathbb{Z}_K=\prod_{i=1}^r\prod_{j=1}^{t_i}\prod_{s=1}^{s_{ij}}\mathfrak{p}_{ijs}^{e_{ij}}.$$
		where $e_{ij}$ is the smallest positive integer satisfying $e_{ij}\la_{ij}\in \mathbb{Z}$ and $f_{ijs}=\deg(\phi_i)\times \deg(\psi_{ijs})$ is the residue degree of $\mathfrak{p}_{ijs}$ over $p$ for every $(i,j,s)$.
	\end{enumerate}
\end{thm}
\begin{cor}\label{cor1}
	Under the hypothesis  above Theorem $\ref{thm4}$, if for every $i=1,\dots,r,\,l_i=1\text{ or }N_{\phi_i}^-(F)=S_i$ has a single side of height $1$, then $\nu_p((\mathbb{Z}_K:\mathbb{Z}[\alpha]))=0$.
\end{cor}
An alternative proof of the theorem of index of Ore is given in \cite[Theorem 1.7 and Theorem 1.9]{EMN}.  In \cite{GMN}, Guardia, Montes, and  Nart introduced  the notion of $\phi$-admissible expansion used in order to treat some special cases when the $\phi$-adic expansion is not obvious. Let
\begin{equation}
\label{eq1}
F(x)=\sum_{i=0}^nA_i'(x)\phi(x)^i,\quad A_i'(x)\in \mathbb{Z}[x],
\end{equation}
be a $\phi$-expansion of $F(x)$, not necessarily deg$(A_i')$ less than deg$(\phi)$. Take $u_i'=\nu_p(A_i'(x))$, for all $i=0,\dots, n$, and let $N'$ be the lower boundary convex envelope of the set of points $\{(i,u_i')\,\mid\,0\leq i\leq n,\,u_i'\neq\infty\}$. To any $i=0,\dots,n$, we attach the residue coefficient as follows:
$$
c_i'=\left\{
\begin{array}{ll}
0,&\text{if }(i,u_i')\text{ lies above }N',\\
\left(\frac{A_i'(x)}{p^{u_i'}}\right)\mod{(p,\phi(x))},& \text{if }(i,u_i')\text{ lies on }N'.
\end{array}
\right.
$$
Likewise, for any side $S$ of $N'$, we can define the residual polynomial attached to $S$ and denoted $R_{\lambda}'(F)(y)\,($similar to the residual polynomial $R_{\lambda}(F)(y)$ from the $\phi$-adic expansion$)$. We say that the $\phi$-expansion $(\ref{eq1})$ is admissible if $c_i'\neq 0$ for each abscissa $i$ of a vertex of $N'$. For more details, we refer to \cite{GMN}.
\begin{lem}\label{NP}$($\cite[Lemma $1.12$]{GMN}$)$\\
	If a $\phi$-expansion of $F(x)$ is admissible, then $N'=N_{\phi}^-(F)$ and $c_i'=c_i$. In particular, for any side $S$ of $N'$ we have $R_{\lambda}'(F)(y)=R_{\lambda}(F)(y)$ up to multiply by a nonzero coefficient of $\mathbb{F}_{\phi}$.
\end{lem}
\section{Proofs of Main Results}

\textbf{Proof of Theorem \ref{thm1}} The proof of this theorem can be concluded by Dedekind's criterion. But as the other results are based on Newton polygons, let us use theorem of index with "if and only if" as it is given in  \cite[Theorem 4.18]{GMN}, which says that a necessary and sufficient condition  to have $\nu_p(\mathbb{Z}_K:\mathbb{Z}[\al])=0$ is that  $ind_1(F)=0$, where $ind_1(F)$ is the index obtained by Ore's index in  Theorem \ref{thm4}. Since $\Delta(F)=\mp 60^{60}\cdot m^{59}$, thanks to the known formula $\nu_p(\Delta(F))=\nu_p(d_K)+2\nu_p(ind(\alpha))$, $\mathbb{Z}[\alpha]$ is the ring of integers of $K$ if and only if $p$ does not divide $(\mathbb{Z}_K:\mathbb{Z}[\alpha])$ for every rational prime $p$ dividing $30\cdot
m$. Let $p$ be a rational prime dividing $m$, then $F(x)\equiv x^{60}\md{p}$. Let $\phi=x$. As $m$ is square free integer, then $\nu_p(m)=1$, and so $N_{\phi}(F)=S$ has a single side  of height $1$. Thus $R_{\lambda}(F)(y)$ is irreducible  over $\mathbb{F}_{\phi}$ as it is of degree $1$. By Corollary $\ref{cor1}$, we get $\nu_p((\mathbb{Z}_K:\mathbb{Z}[\alpha]))=0$; $p$ does not divide $(\mathbb{Z}_K:\mathbb{Z}[\alpha])$. 
For $p=2$ and $2$ does not divide $m$, we have $F(x)\equiv x^{60}-1\equiv (x^{15}-1)^{4}\md{2}$. Let $\phi\in\mathbb{Z}[x]$ be a monic polynomial, whose reduction modulo $2$ is an irreducible factor of $\overline{F(x)}$. Then $\overline{\phi}$ divides $x^{15}-1$ in $\mathbb{F}_2[x]$. Let $x^{15}-1=Q(x)\phi(x)+T(x)$, with $(Q,T)\in \mathbb{Z}[x]^2$ and $\nu_2(T)\ge 1$. Since $x^{15}-1$ is separable over $\mathbb{F}_2$, $\overline{\phi}$ does not $\overline{Q}$.
Let   $F(x)=(x^{15}-1)^4+4(x^{15}-1)^3+6(x^{15}-1)^2+4(x^{15}-1)+1-m=Q(x)^4\phi^4+4Q(x)^3\phi^3+6Q(x)^2\phi^2+4Q(x)\phi+r_0+1-m$, where   $r_0$ is the remainder upon the Euclidean division of $F(x)-(Q(x)^4\phi^4+4Q(x)^3\phi^3+6Q(x)^2\phi^2+4Q(x)\phi+1-m)=2T(x)K(x)$ for some $K\in \mathbb{Z}[x]$,  we conclude that $\nu_2(r_0)\ge 2$. Since $\overline{\phi}$ does not $\overline{Q}$, the previous $\phi$-expansion is admissible, and by Lemma \ref{NP},  {$ind_\phi(F)=0$} if and only if $\nu_2(1-m)=1$; $m\not\equiv 1\md{4}$.\\

 Similarly, for $p=3$ and $3$ does not divide $m$, 
 we have $F(x)\equiv x^{60}-m\equiv (x^{20}-m)^{3}\md{3}$. Let $\phi\in\mathbb{Z}[x]$ be a monic polynomial, whose reduction modulo $3$ is an irreducible factor of $\overline{F(x)}$. Then $\overline{\phi}$ divides $x^{20}-m$ in $\mathbb{F}_3[x]$. 
 Let $x^{20}-m=Q(x)\phi(x)+T(x)$, with $(Q,T)\in \mathbb{Z}[x]^2$ and $\nu_3(T)\ge 1$. Since $x^{20}-m$ is separable over $\mathbb{F}_3$, $\overline{\phi}$ does not $\overline{Q}$.
Let   $F(x)=(x^{20}-m)^3+3m(x^{15}-m)^2+3m^2(x^{15}-m)+m^3-m=Q(x)^3\phi^3+3mQ(x)^2\phi^2+3m^2Q(x)\phi+r_0+m^3-m$, where   $r_0$ is the remainder upon the Euclidean division of $F(x)-(Q(x)^3\phi^3+3mQ(x)^2\phi^2+3m^2Q(x)\phi+m^3-m)=3T(x)K(x)$ for some $K\in \mathbb{Z}[x]$,  we conclude that $\nu_3(r_0)\ge 2$. Since $\overline{\phi}$ does not $\overline{Q}$, the previous $\phi$-expansion is admissible, and by Lemma \ref{NP},   {$ind_\phi(F)=0$} if and only if $\nu_3(m^2-1)=1$; $m^2\not\equiv 1\md{9}$.  \\
For $p=5$ and $5$ does not divide $m$,  let $\phi\in\mathbb{Z}[x]$ be a monic polynomial, whose reduction modulo $5$ is an irreducible factor of $\overline{F(x)}$. Then $\overline{\phi}$ divides $x^{12}-m$ in $\mathbb{F}_5[x]$. 
 Let $x^{12}-m=Q(x)\phi(x)+T(x)$, with $(Q,T)\in \mathbb{Z}[x]^2$ and $\nu_5(T)\ge 1$. Since $x^{12}-m$ is separable over $\mathbb{F}_5$, $\overline{\phi}$ does not $\overline{Q}$.
Let   $F(x)=(x^{12}-m)^5+5m(x^{12}-m)^4+10m^2(x^{12}-m)^3+10m^3(x^{12}-m)^2+5m^4(x^{12}-m)+m^3-m=Q(x)^5\phi^5+5mQ(x)^4\phi^4+10m^2Q(x)^3\phi^3+10m^4Q(x)^2\phi^2+5m^4Q(x)\phi+r_0+m^5-m$, where   $r_0$ is the remainder upon the Euclidean division of $F(x)-(Q(x)^5\phi^5+5mQ(x)^4\phi^4+10m^2Q(x)^3\phi^3+10m^4Q(x)^2\phi^2+5m^4Q(x)\phi+m^5-m)=3T(x)K(x)$ for some $K\in \mathbb{Z}[x]$,  we conclude that $\nu_5(r_0)\ge 2$. Since $\overline{\phi}$ does not $\overline{Q}$, the previous $\phi$-expansion is admissible, and by Lemma \ref{NP},   {$ind_\phi(F)=0$} if and only if $\nu_5(m^4-1)=1$; $m\not\equiv \pm 1, \pm 7 \md{25}$. \\
 \smallskip
 
The index of a field $K$ is defined by $i(K)=gcd\{(\mathbb{Z}_K:\mathbb{Z}[\theta])\mid K=\mathbb{Q}(\theta) \mbox{ and } \theta\in \mathbb{Z}_K \}$. A rational prime $p$ dividing $i(K)$ is called a prime common index divisor of $K$. If $\mathbb{Z}_K$ has a power integral basis, then $i(K)=1$. Therefore a field having a prime common index divisor is not monogenic. For the proof of Theorem $\ref{thm2}$, we need the following lemma and its proof is an immediate consequence of  Dedekind's theorem.
\begin{lem}
	\label{index}
	Let $p$ be a rational prime integer and $K$ be a number field. For every positive integer $f$, let $\mathcal{P}_f$ be the number of distinct prime ideals of $\mathbb{Z}_K$ lying above $p$ with residue degree $f$ and $\mathcal{N}_f$ be the number of monic irreducible polynomials of $\mathbb{F}_p[x]$ of degree $f$. If $\mathcal{P}_f>\mathcal{N}_f$ for some positive integer $f$, then $p$ is a prime common index divisor of $K$.\\
\end{lem}
To apply the last Lemma one has to know the number $\mathcal{N}_f(p)$ of monic irreducible polynomials over $\mathbb{F}_p$ of degree $f$. This number was found by Gauss, which it is given by the following proposition:
\begin{prop}$($\cite[ Chapter 4, Proposition 4.35]{W}$)$\label{pro1}\\
	For every prime $p$ and $f\geq 1$ one has
	$$\mathcal{N}_f(p)=\frac{1}{f}\sum_{d\mid f}\mu(d)p^{f/d},$$
	where $\mu(d)$ is the familiar M\"obius function.
\end{prop}
\textbf{Proof of Theorem \ref{thm2}} In every case, let us show that $i(K)>1$, and so $K$ is not monogenic.
\begin{enumerate}
	\item 
	$m\equiv 1 \md{4}$. Then $\overline{F(x)}=\overline{(x^{15}-1)}^4=\overline{(x-1)(x^2+x+1)U(x)}^4$ in $\mathbb{F}_2[x]$. 
	Let  $\phi(x)=x^2+x+1$ and  $\nu=\nu_2(1-m)$. Since   
	$F(x)=\dots+ (-48165-42465x)\phi(x)^4+(3610+6840x)\phi(x)^3-570x\phi(x)^2+(-20+20x)\phi(x)+1-m$,  	
	if $\nu\geq 4$, then $N_{\phi}^-(F)$ has  $3$ sides joining $(0,\nu),\, (1,2),\, (2,1)$, and $(4,0)$. Thus every side of $N_{\phi}^-(F)$  has degree $1$ (see $Figure\ 2,\, \nu\geq 4$). Thus by Theorem \ref{thm4}, $\phi$ provides $\textcolor{purple}{3}$ prime ideals of $\mathbb{Z}_K$ of residue degree $2$ each one. As there are only one monic irreducible polynomial of degree $2$ in $\mathbb{F}_2[x]$ namely $\overline{\phi(x)}$, by Lemma \ref{index}, $2$ is a common index divisor of $K$, and so $K$ is not monogenic. 
	If $\nu=2$, $N_{\phi}^-(F)=S$ is one sided of degree $2$ such that $R_{\lambda}(F)(y)=(x+1)y^2+xy+1=((x+1)y+1)(y+1)\in\mathbb{F}_{\phi}[y]$ (see $Figure\ 2,\, \nu=2$). So, by Theorem $\ref{thm4}$, $\phi$ provides $2$ prime ideals of $\mathbb{Z}_K$ lying above $2$ with residue degree $2$ each one, and so there are at least  $2$ prime ideals of $\mathbb{Z}_K$ lying above $2$ with residue degree $2$ each one. As there are only one monic irreducible polynomial of degree $2$ in $\mathbb{F}_2[x]$ namely $\overline{\phi(x)}$, by Lemma \ref{index}, $2$ is a common index divisor of $K$, and so $K$ is not monogenic.
	If $\nu=3$, then $N_{\phi}^-(F)=S_1+S_2$ has two sides with degrees $d(S_1)=2$ and $d(S_2)=1$ such that $R_{\lambda_{1}}(F)(y)=xy^2+(1+x)y+1=(y+1)(xy+1)\in \mathbb{F}_{\phi}[y]$ (see $Figure\ 2,\, \nu=3$). Thus, by Theorem \ref{thm4}, $\phi$ provides $3$ prime ideals of $\mathbb{Z}_K$ lying above $2$ with residue degree $2$ each one. As there are only one monic irreducible polynomial of degree $2$ in $\mathbb{F}_2[x]$, $2$ is a common index divisor of $K$, and so $K$ is not monogenic.
	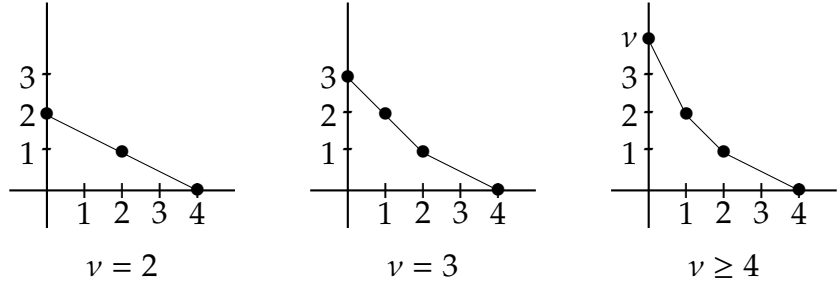
\begin{figure}[htbp] 
		\centering
		\begin{tikzpicture}[x=1cm,y=1cm]
		\draw[thick] (-0.5,0) -- (2.5,0);
		\draw[thick] (0,-0.5) -- (0,2.5);
		\draw (0,0.5) node {-};
		\draw (0,0.5) node[left]{$1$};
		\draw (0,1) node {-};
		\draw (0,1) node[left]{$2$};
		\draw (0,1.5) node {-};
		\draw (0,1.5) node[left]{$3$};
		\draw (0.5,0) node {$\shortmid$};
		\draw (0.5,0) node[below] {$1$};
		\draw (1,0) node {$\shortmid$};
		\draw (1,0) node[below] {$2$};
		\draw (1.5,0) node {$\shortmid$};
		\draw (1.5,0) node[below] {$3$};
		\draw (2,0) node {$\shortmid$};
		\draw (2,0) node[below] {$4$};
		\draw (2,0) node {$\bullet$};
		\draw (0,1) node {$\bullet$};
		\draw(2,0)--(0,1);
		\draw (1,0.5)node{$\bullet$};
		\draw (1,-1)node{$\nu=2$};
		\draw[thick] (3.5,0) -- (6.5,0);
		\draw[thick] (4,-0.5) -- (4,2.5);
		\draw (4,0.5) node {-};
		\draw (4,0.5) node[left]{$1$};
		\draw (4,1) node {-};
		\draw (4,1) node[left]{$2$};
		\draw (4,1.5) node {-};
		\draw (4,1.5) node[left]{$3$};
		\draw (4.5,0) node {$\shortmid$};
		\draw (4.5,0) node[below] {$1$};
		\draw (5,0) node {$\shortmid$};
		\draw (5,0) node[below] {$2$};
		\draw (5.5,0) node {$\shortmid$};
		\draw (5.5,0) node[below] {$3$};
		\draw (6,0) node {$\shortmid$};
		\draw (6,0) node[below] {$4$};
		\draw (6,0) node {$\bullet$};
		\draw (4,1.5) node {$\bullet$};
		\draw (5,0.5)node {$\bullet$};
		\draw (4.5,1) node {$\bullet$};
		\draw(4,1.5)--(5,0.5);
		\draw (5,0.5)-- (6,0);
		\draw (5,-1)node{$\nu=3$};
		\draw[thick] (7.5,0) -- (10.5,0);
		\draw[thick] (8,-0.5) -- (8,2.5);
		\draw (8,0.5) node {-};
		\draw (8,0.5) node[left]{$1$};
		\draw (8,1) node {-};
		\draw (8,1) node[left]{$2$};
		\draw (8,1.5) node {-};
		\draw (8,1.5) node[left]{$3$};
		\draw (8.5,0) node {$\shortmid$};
		\draw (8.5,0) node[below] {$1$};
		\draw (9,0) node {$\shortmid$};
		\draw (9,0) node[below] {$2$};
		\draw (9.5,0) node {$\shortmid$};
		\draw (9.5,0) node[below] {$3$};
		\draw (10,0) node {$\shortmid$};
		\draw (10,0) node[below] {$4$};
		\draw (8,2) node {$\bullet$};
		\draw (10,0) node {$\bullet$};
		\draw (8.5,1)node{$\bullet$};
		\draw(8.5,1)--(8,2);
		\draw (9,0.5)node{$\bullet$};
		\draw (8,2)node[left]{$\nu$};
		\draw (8.5,1)--(9,0.5);;
		\draw (9,0.5)--(10,0);
		\draw (9,-1)node{$\nu\geq 4$};
		\end{tikzpicture}
		\caption{\large  The principal $\phi$-Newton polygon of $F(x)$.}
	\end{figure}
	\item 
	$m\equiv \mp 1\md{9}$.\\
	For 
	$m\equiv  1\md{9}$,  since $\overline{F(x)}=\overline{\phi_1(x)\phi_2(x)U(x)}^3$  in {$\mathbb{F}_3[x]$}, with  $\phi_1(x)=x-1$, $\phi_2(x)=x+1$ and $\overline{\phi_k}$ does not divide $\overline{U(x)}$ in $\mathbb{F}_3[x]$ for every $k=1,2$.
	Consider the following expansions: $F(x)=\dots+34220\phi_1(x)^3+1770\phi_1(x)^2+60\phi_1(x)+1-m$ and  $F(x)=\dots-34220\phi_2(x)^3+1770\phi_2(x)^2-60\phi_2(x)+1-m$. We conclude that if $m\equiv 1\md{9}$, then	$N_{\phi_k}^-(F)=S_{k1}+S_{k2}$ has $2$ sides joining $(0,V)$, $(1,1)$, and $(3,0)$ with $V=\nu_3(1-m)\geq  2$. Thus the degree of each side is $1$. Therefore, $\phi_k$ provides $2$ prime ideals of $\mathbb{Z}_K$ lying above $3$ with residue degree $1$ each one. Applying this for every $k=1,2$, we conclude that there are $4$ prime ideals of $\mathbb{Z}_K$ lying above $3$ of residue degree $1$ each one. As there are only $3$ monic irreducible polynomials of degree $1$ in $\mathbb{F}_3[x]$, $3$ is a common index divisor of $K$ and so $K$ is not monogenic.\\
	For $m\equiv -1\md{9}$,  we have  since $\overline{F(x)}=\overline{(x^{20}+1)}^3=\overline{\phi_1(x)\phi_2(x)U(x)}^3\md{3}$ with $\phi_1(x)=x^2+x-1$, $\phi_2(x)=x^2-x-1$, and $\overline{\phi_k}$ does not divide $\overline{U(x)}$ in $\mathbb{F}_3[x]$ for every $k=1,2$.
	Consider the following expansions: $F(x)=\dots+a_3(x)\phi_1(x)^3+a_2(x)\phi_1(x)^2+a_1(x)\phi_1(x)+a_0(x)$ and\\ $F(x)=\dots+b_3(x)\phi_2(x)^3+b_2(x)\phi_2(x)^2+b_1(x)\phi_2(x)+b_0(x)$, where $a_1(x)=16175489617620-25052342327220x$, $a_0(x)=-1548008755920x+956722026041-m$, $b_1(x)=16175489617620+25052342327220x$, and $b_0(x)=1548008755920x+956722026041-m$. Since $\nu_3(a_3(x)b_3(x))=0$,  $\nu_3(a_2(x))\geq 1$, {$\nu_3(b_2(x))\geq 1$},  {$\nu_3(a_1(x))=\nu_3(b_1(x))=1$}, we conclude that if $m\equiv -1\md{9}$, then $\nu_3(a_0(x))\geq 2$ and $\nu_3(b_0(x))\geq 2$, and so	$N_{\phi_k}^-(F)=S_{k1}+S_{k2}$ has $2$ sides joining $(0,V_k)$, $(1,1)$, and $(3,0)$ with {$V_1=\nu_3(a_0(x))\geq  2$ and $V_2=\nu_3(b_0(x))\geq  2$}. Thus the degree of each side is $1$. Therefore, $\phi_k$ provides $2$ prime ideals of $\mathbb{Z}_K$ lying above $3$ with residue degree $2$ each one.  The total is then, there are $4$ prime ideals of $\mathbb{Z}_K$ lying above $3$ of residue degree $2$ each one. As the are only $3$ monic irreducible polynomials of degree $2$ in $\mathbb{F}_3[x]$ namely $\overline{\phi_1(x)},~\overline{\phi_2(x)}\text{ and }x^2+1$, $3$ is a common index divisor, and so $K$ is not monogenic.
	\item 
	For $m\equiv \mp 1\md{25}$,\\
	{If $m\equiv 1\md{25}$, we have} ${\overline{F(x)}={(x^{12}-1)}^5}={\prod_{k=1}^4\overline{\phi_k(x)U(x)}^5}$ in $\mathbb{F}_5[x]$, with $\phi_k(x)=x-k$ for every $k=1,\dots,4$ and $\overline{\phi_k(x)}$ does not divide $\overline{U(x)}$ in $\mathbb{F}_5[x]$.
Similarly, by considering the $\phi_k$-expansion of $F(x)$, if $\nu_5(m-1)\geq 2$; {$m\equiv 1\md{25}$}, then  $N_{\phi_k}^-(F)=S_{k1}+S_{k2}$ has $2$ sides joining $(0,V_{k}),~(1,1)\text{ and }(5,0)$, with $V_k\geq 2$. Thus  each side of   $N_{\phi_k}^-(F)$ is of degree $1$. Therefore $\phi_k$ provides $2$ prime ideals of $\mathbb{Z}_K$ lying above $5$ with residue degree $1$ each one. Apply this  for every $k=1,\dots,4$, we conclude that there are $8$ prime ideals of $\mathbb{Z}_K$ lying above $5$ of residue degree $1$ each one. As there are only $5$ monic irreducible polynomials of degree $1$ in $\mathbb{F}_5[x]$, by Lemma $\ref{index}$, $5$ is a common index divisor, and so $K$ is not monogenic.\\
	If $m\equiv -1\md{25}$, then $\overline{F(x)}={x^{12}+1}^5=\prod_{k=1}^6\phi_k(x)U(x)\md{5}$, with $\phi_1(x)=x^2+2,~\phi_2(x)=x^2+3,~\phi_3(x)=x^2+x+2,~\phi_4(x)=x^2+2x+3,~\phi_5(x)=x^2+3x+3,~\phi_6(x)=x^2+4x+2\text{ and }\overline{\phi_k}\nmid \overline{U(x)}$ for every $k=1,\dots,6$.
	Fix {$k=1,\dots,6$} and consider the $\phi_k$-expansion of $F(x)$. If $\nu_5(m+1)\geq 2$, then $N_{\phi_k}^-(F)=S_{k1}+S_{k2}$ has $2$ sides joining $(0,V_k),~(1,1)\text{ and }(5,0)$ with $V_k\geq 2$. Thus every side of $N_{\phi_k}^-(F)$ is of  degree $1$. It follows by Theorem  \ref{thm4} that  $\phi_k$ provides $2$ prime ideals of $\mathbb{Z}_K$ lying above $5$ with residue degree $2$ each one.
	Applying this for every  $k=1,\dots, 6$, we conclude that there are $12$ prime ideals of $\mathbb{Z}_K$ lying above $5$ of residue degree $2$ each one. By Proposition $\ref{pro1}$ the are only $10$ monic irreducible polynomials of degree $2$ in $\mathbb{F}_5[x]$, by Lemma $\ref{index}$, $5$ is a common index divisor, and so $K$ is not monogenic.
\end{enumerate}
\begin{rem}
Let $F(x)=x^n-m\in\mathbb{Z}[x]$ be an irreducible  polynomial over $\mathbb{Q}$ and $K=\mathbb{Q}(\alpha)$ with $\alpha$ a complex root of $F(x)$. Let $p$ be a prime integer dividing $n$ and does not divide $m$, and let $r=\nu_p(n)$.
In \cite{Gs}, Gassert claimed that $N_{\phi}^-(F)$ is the convex envelope of the set of points $\{(0,\nu_p(m^p-m))\}\cup\{(k,\nu_p(\binom{p^r}{k})),\, k=1,\dots,r\}$. 
The following example shows that this claim is not correct.
$F(x)=x^{60}-m$ with $m\neq \pm 1$ a square free integer such that $m\equiv -1\md{27}$. Then for $p=3$ and $\phi=x^2+x-1$, we have $F(x)\dots+a_3(x)\phi_1(x)^3+a_2(x)\phi_1(x)^2+a_1(x)\phi_1(x)+a_0(x)$, with $a_1(x)=16175489617620-25052342327220x$, $a_0(x)=-1548008755920x+956722026041-m$.
 As $\nu_3(a_0(x))=2$, then $N_{\phi}^-(F)$ is the convex envelope of the set of points $\{(0,2), (1,1), (3,0)\}$ contrary to the claim, which says that it will be  the convex envelope of the set of points $\{(0,V), (1,1), (3,0)\}$ with $V\ge 3$. 
\end{rem}
\textbf{Proof of Theorem \ref{cor}} As gcd$(u,30)=1$, let $(x,y)\in\mathbb{Z}^2$ be  the unique solution of  $ux-60y=1$ with $0\le y<u$ and let $\theta =\frac{\alpha^x}{a^y}$. Then $\theta^{60}=\frac{\alpha^{60 x}}{a^{60y}}=a^{ux-60 y}=a$. Since $g(x)=x^{60}-a\in \mathbb{Z}[x]$ is an Eisenstein polynomial, $g(x)$ is irreducible over $\mathbb{Q}$. As $\theta\in K$ and $[K:\mathbb{Q}]=\deg(g)$, we conclude that $K=\mathbb{Q}(\theta)$. Therefore, $K$ is generated by a root of the polynomial $g(x)=x^{60}-a$ with $a\neq \mp 1$  a square free integer.  The proof is therefore  an application of  Theorem \ref{thm1} and  Theorem \ref{thm2}.
\begin{flushright}
	$\square$
\end{flushright}
In order to illustrate the efficiency of our results, we finalize the paper by the following numerical  examples.
\begin{example}
	Let $F(x)\in \mathbb{Z}[x]$ be a monic irreducible polynomial and $K$ the number field generated by a complex root of $F(x)$.
	\begin{enumerate}
		\item If $F(x)=x^{60}-67$, then $F(x)$ is irreducible because it is $67$-Eisenstein. Since $m\equiv 3\md{4}$, $m\equiv 4 \md{9}$ and $m\equiv 17 \md{25}$, by Theorem $\ref{thm1}$ $K$ is monogenic.
		\item If $F(x)=x^{60}-302$, then $F(x)$ is irreducible because it is $2$-Eisenstein. Since $m\equiv 2\md{4}$, $m\equiv 5 \md{9}$ and $m\equiv 6 \md{25}$, by Theorem $\ref{thm1}$ $K$ is monogenic.
		\item If $F(x)=x^{60}-106$, then $F(x)$ is irreducible because it is $2$-Eisenstein. Since $m\equiv 1\md{5}$, by Theorem $\ref{thm2}$ $K$ is not monogenic.
		\item If $F(x)=x^{60}-226$, then $F(x)$ is irreducible because it is $2$-Eisenstein. Since $m\equiv 1\md{9}$, by Theorem $\ref{thm2}$ $K$ is not monogenic.
		\item If $F(x)=(x-5)^{60}-70^{13}$, then $F(x)\equiv x^{60}\md{5}$. As $70\equiv 2\md{4},~70\equiv 7\md{9}\text{ and }70\equiv -5\md{25}$, by Theorem $\ref{cor}$, $K$ is monogenic.
		\item If $F(x)=(x-4)^{60}-26^{31}$, then $F(x)\equiv x^{60}\md{2}$. As $26\equiv 1\md{25}$, by Theorem $\ref{cor}$, $K$ is not monogenic.
	\end{enumerate}
\end{example}
\begin{rem}
In all calculations of $\phi$-expansions, we used Maple 12.
\end{rem}

\end{document}